%% file: omega-n-y.tex
\magnification 1050
\input itenn.tex

\input option_keys.tex
\ifx\montrerlabels\oui
\input montrerlabels.tex
\fi
\newcount\numchapi\numchapi=0

\optionparag=1
\def\paradouze{oui}
\anglais
\dimstand\hautspages{G. Tenenbaum}{On local laws for the number of small prime factors}
{\leftskip-10mm\obeylines
\it Ramanujan Journal 
\bf 51 \rm(2020), 153--161.\par }
\titrecentre{On local laws for the number of small prime factors\anote{*}{\hskip-4mm We include here the correction of a misprint appearing in the statement of Theorem 1.5 of the published version.\hfill}}
\bigskip\medskip
\centerline{GŽrald Tenenbaum}
\bigskip\bigskip
{\eightpoint\leftskip1cm\rightskip1cm
\noi{\bf Abstract.} Investigating  a question of Alladi, we describe the local distribution of small prime factors of integers, with emphasis on the  transition phase occurring for certain values of the parameters.
\PAR
\medskip\noi
\par
\smallskip 
\noi \bf 2010 Mathematics Subject Classification: \rm   11N25, 11N37, 11N60.\par }
\par 
\bigskip\bigskip
\paraun{Introduction and statement of results}
Let $\nu(n,y):=\sum_{p|n,\,p\leqslant y}1$ denote the number of prime factors not exceeding $y$, counted without multiplicity, of a natural integer $n$. K. Alladi proposed to estimate 
$$N_k(x,y):=\sum_{\di{n\leqslant x}{\nu(n,y)=k}}1\qquad (x\geqslant y\geqslant 2, \, k\ll\log_2y)$$
with the purpose of describing the phase transition between the case of small $y$, when we expect
$$N_k(x,y)\asymp {x(\log_2y)^k\over k!\log y}\eqdef{faNksmall}$$
and the case of large $y$, when the parameter $k$, in the above estimate, should be replaced by $k-1$.\note{Here and in the sequel we let $\log_j$ denote the $j$-th iterated logarithm.}
Of course, as stated, this phenomenon can only be brought to light when $r:=k/\log_2y$ is small. 
\par 
As we shall see, this question is quite appealing and a complete answer still eludes us. We propose here, as a benchmark, a first, concise study resting principally on results derived through the saddle-point method. In parallel to the present work, Alladi and Molnar \citer{AM18} have tackled the problem  using mainly Buchstab's iteration method.\par 
We  observe incidentally that the situation changes in nature when $k/\log_2y$ becomes large. For instance, we have from corollary 1 in \citer{HT89},
$$\pi_k(x):=N_k(x,x)\sim {F(r)\over \Gamma(r+1)}{\e^{-kh/2}x(\log_2x)^{k-1}\over (k-1)!\log x}=o\Big({x(\log_2x)^{k-1}\over (k-1)!\log x}\Big)$$
if $r:=k/\log_2x\to\infty$, $k=o\big((\log_2x)^2\big)$, with 
$$F(z):=\prod_{p}\Big(1+{z\over p-1}\Big)\Big(1-{1\over p}\Big)^z,\quad h=\bigg({\log (2+\e^\gamma y\log y)\over \log_2x}\bigg)^2.\eqdef{Fz}$$\par 
In the case of relatively small values of $y$, an asymptotic formula with remainder easily follows from estimates established in \citer{ET89}. Here and in the sequel, $\gamma$ denotes Euler's constant.
\Propt{thmk}{There exists an absolute constant $c>0$ such that,  uniformly under the conditions $3\leqslant y\leqslant x^{c(\log_3x)/\log_2x}$, $r:=k/\log_2y\ll1$, we have
$$N_k(x,y)=F(r)\e^{\gamma(r-1)}{x(\log_2y)^k\over k!\log y}\Big\{1+O\Big({k\over (\log_2 y)^2}\Big)\Big\}\cdot\eqdef{f1bis}$$}
\medskip
Next, we investigate the case of larger values of $y$. The following result is an easy consequence of a special case of corollary 2.4 in \citer{Te17}. We write
$$S_z(x,y):=\sum_{n\leqslant x}z^{\nu(n,y)}\qquad (x\geqslant y\geqslant 1,\,z\in\CC).$$
\par 
\Propt{univ}{Let $\kappa>0$. Uniformly for $3\leqslant y\leqslant x$, $1\leqslant k\leqslant (\log_2y)/\kappa$, we have
$$N_k(x,y)\ll x{(\log_2y)^k\over k!\,\log y}\cdot\eqdef{majuniv}$$
Moreover, for any fixed $\kappa\in]0,1[$ and uniformly for $3\leqslant y\leqslant x$, $\kappa\leqslant r:=k/\log_2 y\leqslant 1/\kappa$, we have
$$N_k(x,y)=S_r(x,y){(\log_2y)^k\over k!\e^k}\Big\{1+O\Big({1\over \sqrt{\log_2y}}\Big)\Big\}.\eqdef{Mk/Sr}$$}
\par 
In order to exploit \eqref{Mk/Sr} we need to estimate $S_r(x,y)$. This depends on the quantity
$$\sigma_r(u):=\int_0^{u-1}\omega(u-t)\varrho_r(t)\d t+\varrho_r(u)\quad(u\geqslant 1),$$
where $\varrho_r$ is the $r$-th fractional convolution power of the Dickman function $\varrho=\varrho_1$---it satisfies $\varrho_r(u)=\varrho(u)r^{u+o(u)}$ as $u\to\infty$, see \citer{TW03} for further details---and $\omega$ is Buchstab's function---see, \eg, \citer{Te15}, \S III.6.3. Appealing to the estimate $\omega(t)=\e^\gamma+O(t^{-t})$, proved for instance in \citer{Te15}, th. III.6.8, we easily obtain that for fixed $r>0$,
$$\sigma_r(u)>0\quad(u>1),\qquad \sigma_r(u)=\e^{\gamma(r-1)}+O(u^{-u})\qquad (u\to\infty).\eqdef{sigr}$$\par 
For fixed $\varepsilon>0$, we define the
domain $$\exp\big\{(\log_2x)^{5/3+\varepsilon}\big\}\leqslant y\leqslant x.\leqno{(H_\varepsilon)}$$
\par 
\Propt{H+sqrt}{Let $\varepsilon>0$. Uniformly for $(x,y)\in(H_\varepsilon)$,  $3\leqslant y\leqslant \sqrt{x}$, $r:=k/\log_2y\ll1$, and with $u:=(\log x)/\log y$, we have
$$S_r(x,y)=F(r)\sigma_r(u)x(\log y)^{r-1}\Big\{1+O\Big({1\over \log y}\Big)\Big\}.\eqdef{faSr}$$
Consequently, for any fixed $\kappa\in]0,1[$ and uniformly for $\sqrt{x}\geqslant y\geqslant 3,$ $\kappa\leqslant r\leqslant 1/\kappa,$ we have
$$N_k(x,y)={F(r)\sigma_r(u)x(\log_2y)^k\over k!\log y}\Big\{1+O\Big({1\over\sqrt{\log_2y}}\Big)\Big\}.\eqdef{Mk/Sr2}$$
}
 \par 
 Note that, in view of \eqref{sigr}, formulae \eqref{Mk/Sr2} and \eqref{f1bis} coincide for large $u$.
 \par 
 At this point, we are left with two cases : \par 
 {\parindent10mm (a) $x^{c(\log_3x)/\log_2x}<y\leqslant \sqrt{x}$ and $k\leqslant \kappa\log_2y$ with fixed, arbitrary small $\kappa>0$; \par  (b) $\sqrt{x}<y\leqslant x$, $1\leqslant k\leqslant (\log_2y)/\kappa$.\par} 
 We first show that  behaviour \eqref{faNksmall} holds throughout range (a), and somewhat beyond.
 \Propt{Mk-ypt}{Let $c\in]0,1[$ and $\kappa>0$ be fixed. Then, estimate \eqref{faNksmall} holds uniformly for $2\leqslant y\leqslant x^{1-c}$, $1\leqslant k\leqslant (\log_2y)/\kappa$.}
\par \smallskip
From the above, we now know that the desired threshold must occur in range (b). 
Our final, general statement takes into account the case when $w:=(\log x)/\log (3x/y)$ is  large. 
\Propt{th-wgd}{Let $\kappa>0.$ Uniformly for $x\geqslant y\geqslant 3$, $1\leqslant k\leqslant (\log_2y)/\kappa$, we have
$$N_k(x,y)\asymp  {x(\log_2y)^{k-1}\over (k-1)!\log y}+{x(\log_2\min(3x/y,y))^k\over k!\log y}\cdot \eqdef{faMk-2}$$}
It follows in particular from this statement that, when $w\geqslant 3$, the estimate \eqref{faNksmall} holds 
if $1\leqslant k\ll (\log_2y)/\log w $, while we have
$$N_k(x,y)\asymp  {x(\log_2y)^{k-1}\over (k-1)!\log y}\eqdef{faMk-3}$$ when  $1+(\log_2y)(\log_2w)/\log w\leqslant k\ll\log_2y$. In the middle range, formula \eqref{faMk-2} describes the transition.
\goodbreak
\medskip
We have thus determined the true order of magnitude of $N_k(x,y)$ for large $y$ and all $k\ll\log_2y$, and, in a number of cases, provided an asymptotic formula with remainder.
It is a challenging problem, with interesting methodological issues, to obtain asymptotic estimates in the remaining ranges.

\bigskip\goodbreak
\paraun{Proof of \ref{thmk}}
Let $z\in\CC$. By a slight modification of the proof of theorem 02 of \citer{HT88} (\ie appealing to classical upper bounds for the number $\Psi(x,y)$ of $y$-friable integers not exceeding $x$ instead of the simple Rankin bound used in \citer{HT88})  we get, uniformly for $2\leqslant y\leqslant x^{c(\log_3x)/\log_2x}$, $z\ll1$,
$$S_z(x,y)=x\prod_{p\leqslant y}\bigg(1+{z-1\over p}\Big)+O\Big({x\over (\log x)^2}\bigg).$$
Since the main term of \eqref{f1bis} is of order $\gg x/\log y$, we see that $N_k(x,y)/x$ is, to the stated accuracy, approximated by the coefficient of $z^k$ in the product over primes. However, this has been studied in \citer{ET89}. 
 For $k\ll\log_2y$, let us write
$$L:=\log \Big({\log y\over \log (k+1)}\Big), \quad M:=\log_2y-\log \big(1+\log^+(k/L)\big)=\log_2y+O(1), \quad\varrho:=k/M.$$
Then, by corollary 1 of \citer{ET89},
$$N_k(x,y)=x\prod_{p\leqslant y}\Big(1+{\varrho-1\over p}\Big){M^k\over k!\,\e^k}\Big\{1+O\Big({k\over (\log_2 y)^2}\Big)\Big\}\cdot\eqdef{f1}$$
Formula \eqref{f1bis} follows by writing $\varrho=r+O\big(r/\log_2y\big)$ and noticing that the main terms of \eqref{f1bis} and \eqref{f1} coincide to within the stated accuracy.
\par 
\smallskip
\bigskip
\paraun{Proof of \ref{univ}}
The upper bound\eqref{majuniv} follows immediately from corollary 2.4(iii) of \citer{Te17}. The same statement provides \eqref{Mk/Sr} with $E(y):=\sum_{p\leqslant y}1/p$ instead of $\log_2y$ and $k/E(y)$ instead of~$r$. However, shifting $E(y)$ by a bounded amount does not perturb the asymptotic formula as stated. This can be seen in either of two ways. The first is by directly inspecting the proof displayed in \citer{Te17}, where a quantity $\e^{(z-r)E(y)}\{1+O(\vartheta)\}$ gets integrated over the circle $z=r\e^{i\vartheta}$; the second is by using the simple estimate
$$S_r(x,y)\asymp x(\log y)^{r-1}\qquad (r\asymp1),\eqdef{estSr}$$
which readily follows, for instance, from theorem 1.1 of \citer{Te16}. We only outline this approach. Applying \eqref{estSr} with $(1\pm v)r$ in place of $r$ and choosing $v$ optimally,
we obtain, for bounded, positive $v$, 
$${1\over S_r(x,y)}\sum_{\di{n\leqslant x}{|\nu(n,y)-r\log_2y|>vr\log_2y}}r^{\nu(n,y)}\ll(\log y)^{-rQ(v)}$$
with $Q(v)\asymp v^2$ if $0\leqslant v\leqslant 1/2$, and $Q(v)\asymp (1+v)\log (1+v)$ if $v>1/2$. Integrating over~$v$, we obtain
$$\sum_{n\leqslant x}r^{\nu(n,y)}\{\nu(n,y)-r\log_2y\}^2\ll S_r(x,y)\log_2y.\eqdef{TKSr}$$
From this, we see that $S_r(x,y)r^{-k}=\sum_{n\leqslant x}r^{\nu(n,y)-k}$ varies by at most a factor $$1+O(1/ \sqrt{\log_2y})$$ when $r-k/\log_2y\ll1/\log_2y$.
\par 
We note in passing that \eqref{TKSr} also follows from the weighted version of the Tur‡n-Kubilius inequality proved in \citer{BS96}. However, verifying the hypotheses then turns out to be more complicated.
\paraun{Proof of \ref{H+sqrt}}
Using the classical notation $\Phi(x,y)$ for the number of uncancelled integers in the sieve by prime factors $\leqslant y$ and letting $P^+(n)$ denote the largest prime factor of an integer $n$ with the convention that $P^+(1)=1$, we have
$$S_r(x,y)=T_r(x,y)+U_r(x,y)-U_r(x/y,y),\eqdef{decSr}$$
with
$$T_r(x,y):=\sum_{\di{m\leqslant x/y}{P^+(m)\leqslant y}}r^{\nu(m)}\Phi\Big({x\over m},y\Big), \qquad U_r(x,y):=\sum_{\di{m\leqslant x}{P^+(m)\leqslant y}}r^{\nu(m)}.$$
By a result in \citer{TW03}, writing  $u:=(\log x)/\log y$, we have
$$U_r(x,y)=\Big\{1+O\Big({\log (u+1)\over \log y}\Big)\Big\}F(r)\varrho_r(u)x(\log y)^{r-1}\eqdef{faUr}$$
in $(H_\varepsilon)$,
with notation \eqref{Fz}.\par 
In order to evaluate $T_r(x,y)$, we use the simple formula---see \citer{Te15}, corollary III.6.20---, valid in $(H_\varepsilon)$,
$$\Phi(x,y)=\e^{\gamma}{x\omega(u)-y\over \zeta(1,y)}+O\Big({x\varrho(u)\over (\log y)^2}\Big),\eqdef{appPhi}$$
where $\omega$ is Buchstab's function, $\gamma$ is Euler's constant, and $\zeta(s,y):=\prod_{p\leqslant y}(1-1/p^s)^{-1}$. This yields
$$T_r(x,y)={\e^\gamma x\over \zeta(1,y)}A_r(x,y)-{\e^\gamma y\over \zeta(1,y)}U_r(x/y,y)+O\Big({x\over (\log y)^2}B_r(x,y)\Big)$$
where we have put
$$\eqalign{A_r(x,y)&:=\sum_{\di{m\leqslant x/y}{P^+(m)\leqslant y}}{\omega(u-u_m)r^{\nu(m)}\over m},\quad B_r(x,y):=\sum_{\di{m\leqslant x/y}{P^+(m)\leqslant y}}{\varrho(u-u_m)r^{\nu(m)}\over m},\cr}$$
with $u_m:=(\log m)/\log y$.\par 
Appealing to \eqref{faUr}, partial summation yields
$$\eqalign{A_r(x,y)&=\int_{0-}^{u-1}{\omega(u-v)\over y^v}\d U_r(y^v,y)\cr&={U_r(x/y,y)\over x/y}+\int_0^{u-1}U_r(y^v,y)\Big\{{\omega'(u-v)+\omega(u-v)\log y\over y^v}\Big\}\d v\cr
&=F(r)(\log y)^{r}\int_0^{u-1}\varrho_r(v)\omega(u-v)\d v+O\Big((\log y)^{r-1}\Big)\cr}$$
while we have trivially
$$B_r(x,y)\ll(\log y)^r
$$
Assembling our estimates, we get \eqref{faSr}.
\par 
We deduce \eqref{Mk/Sr2} from \eqref{Mk/Sr} and \eqref{faSr} in $(H_\varepsilon)$, and from \eqref{f1bis} and \eqref{sigr} in the complementary domain.
\paraun{Proof of \ref{Mk-ypt}} 
Put $\nu(n):=\nu(n,n)$. Note that, parallel to \eqref{decSr}, we have the decomposition
$$N_k(x,y)=\sum_{\tri{m\leqslant x/y}{\nu(m)=k}{P^+(m)\leqslant y}}\Phi\Big({x\over m},y\Big)+\Theta_k(x,y),\eqdef{decMk}$$
with 
$$\Theta_k(x,y):=\sum_{\tri{x/y<m\leqslant x}{\nu(m)=k}{P^+(m)\leqslant y}}1.$$
In view of \eqref{majuniv}, we only have to establish the lower bound. For this, we  retain the first term on the right-hand side of \eqref{decMk} and appeal to \eqref{appPhi}, using the fact that $\omega(u)\geqslant \dm$ for $u\geqslant 1$. Writing $z:=\ft13y^c$, so that $x/z\geqslant 3y$, we get
$$N_k(x,y)\gg {x\over \log y}\sum_{\di{m\leqslant z}{\nu(m)=k}}{1\over m}\gg {x(\log_2z)^k\over k!\,\log y}\asymp {x(\log_2y)^k\over k!\,\log y}\cdot\eqdef{minMk}$$
\paraun{Proof of \ref{th-wgd}}
We may plainly assume $w$ to be arbitrarily large.
\par 
The case $k=1$ may be dealt with directly on observing that $n$ is counted by $N_1(x,y)$ if, and only if, $n$ is either a prime power $p^j$ with $p\leqslant y$ or is of the form $p^jq$ where $p$ and $q$ are primes and $p^j\leqslant x/y$, $y<q\leqslant x/p^j$. We omit the details and  assume $k\geqslant 2$ in the sequel.\par 
First consider the case $k\ll (\log_2y)/\log w$, so that \eqref{faMk-2} amounts to \eqref{faNksmall}. In view of \eqref{majuniv}, we thus only have to show the lower bound. However, in the indicated range, this is provided by \eqref{minMk} with now $z:=\ft13y^{1/w}$ so that $x/z\geqslant 3y$ and $\log_2z> \log_2y-\log w+\log_2 2\gg k.$
\par 
Next, we embark on proving \eqref{faMk-2} for the larger values of $k$.  From \eqref{decMk}, we  have
$$\eqalign{N_k(x,y)\leqslant R_k(x,y)+\pi_k(x)\cr}$$
with
$$\eqalign{\pi_k(x)\asymp{x\over \log y}{(\log_2y)^{k-1}\over (k-1)!},\quad R_k(x,y)\asymp\sum_{\di{m\leqslant x/y}{\nu(m)=k}}{x\over m\log y}\ll{x(\log_23x/y+c_1)^{k}\over k!\,\log y}.\cr}$$
Here and in the remainder of this proof, $c_j$ $(j\geqslant 1)$ denotes an absolute constant. The expected upper bound follows on noticing that, on the right-hand side of \eqref{faMk-2}, the order of magnitude of the second term can only exceed that of the first if $k\ll\log_23x/y$.
\par 
We establish the lower bound in three steps.
First assume $k\geqslant 1+(\log_2w)(\log_2 y)/\log w$, so that, as stated in \eqref{faMk-3}, the right-hand side of \eqref{faMk-2} is $\asymp \pi_k(x)$. Observe that, if an integer $n\leqslant x$ satisfies $\nu(n)=k$, then either $n$ is counted by $N_k(x,y)$ or $n$ is divisible by a prime  $p\in]y,x]$. Therefore
$${x(\log_2y)^{k-1}\over (k-1)!\log y}\asymp\pi_k(x)\leqslant N_k(x,y)+U_k(x,y),$$
with 
$$\eqalign{U_k(x,y)&:=\sum_{\tri{y<p\leqslant x}{mp\leqslant x}{\nu(m)=k-1}}1\ll\sum_{\di{m\leqslant x/y}{\nu(m)=k-1}}{x\over m\log y}\cr&\ll{x(\log_23x/y+c_2)^{k-1}\over (k-1)!\log y}\ll{x(\log_2y-\log w+c_3)^{k-1}\over (k-1)!\log y}\cdot\cr}\eqdef{majNk}$$
By our assumption on $k$ (and actually $k\geqslant 1+C(\log_2y)/\log w$ with a suitable large $C$  suffices here), we thus have $U_k(x,y)\leqslant \dm\pi_k(x)$ and the required lower bound for $N_k(x,y)$ follows.    
\par\goodbreak
Next, we prove \eqref{faMk-2} for the range $1+C(\log_2 y)/\log w<k\leqslant 1+(\log_2 y)(\log_2w)/\log w$ and under the extra assumption that  $y\leqslant  x/\log x$. We exploit this last hypothesis in the form that  \hbox{$k\leqslant (\log_2w)(\log_2y)/\log w$} implies $k\ll\log_23x/y$: this is obvious if, say, $w\leqslant \sqrt{\log y}$ and otherwise follows from the extra hypothesis since then $\log_23x/y\gg\log_3y$.  Consider \eqref{decMk}. The first term on the right-hand side is plainly
$$\asymp{x(\log_23x/y)^k\over k!\log y}\cdot$$ 
However, as previously noticed, since $\Theta_k(x,y)+U_k(x,y)=\pi_k(x)$, the upper bound \eqref{majNk} implies that $\Theta_k(x,y)\asymp\pi_k(x)$ provided $C$ is large enough. \par 
Finally, we consider the case $x/\log x<y\leqslant x$ and $k> 1+C(\log y)/\log w$ with some sufficiently large constant $C$. Then the right-hand side of \eqref{faMk-2} is $\asymp\pi_k(x)$ because $k\geqslant 2$.
Since  $\Theta_k(x,y)\gg\pi_k(x)$, the required lower bound follows.
\medskip
\noi{\bf Acknowledgements.} The author wishes to express warm thanks to Krishnaswami Alladi for bringing this problem to his attention and for subsequent interesting discussions on the matter.

\bigskip\bigskip
\centerline{\twelvebf References}\bigskip
\eightpoint{
\bibtem{AM18} K. Alladi \& T. Molnar, On the local distribution of the number of small prime factors, preprint.
\bibtem{BS96} A. Bir\'o \& T. Szamuely, A Tur\'an--Kubilius inequality with multiplicative weights,
 {\it Acta Math. Hung. \bf 70} (1996), \numeros 1-2, 39--56.\par 
\bibtem{ET89} P. Erd\H os \& G. Tenenbaum, Sur les densitŽs de certaines suites d'entiers, {\it
Proc. London Math. Soc.} {\rm (3)} {\bf 59} (1989), 417--438.\par 
\bibtem{HT88} R.R. Hall \& G. Tenenbaum, {\it Divisors}, Cambridge tracts in
mathematics, no 90, Cambridge University Press (1988).\par 
\bibtem{HT89} A. Hildebrand \& G. Tenenbaum, On the number of prime factors of an integer, {\it Duke Math.
J. \bf 56} no 3 (1988), 471--501.\par 
\bibtem{Te15} G. Tenenbaum, {\it Introduction to analytic and probabilistic number theory}, 3rd ed., Graduate Studies in Mathematics 163, Amer. Math. Soc. 2015.\par
\bibtem{Te16} G. Tenenbaum, Fonctions multiplicatives, sommes d'exponentielles, et loi des grands nombres, {\it Indag. Math. \bf27} (2016), 590--600.
\bibtem{Te17} G. Tenenbaum, Valeurs moyennes effectives de fonctions multiplicatives complexes, {\it Ramanujan J. \bf44}, \numero3 (2017), 641-701.\par
\bibtem{TW03} G. Tenenbaum \& J. Wu, Moyennes de certaines fonctions
multiplicatives sur les entiers friables, {\it J. reine angew. Math. \bf 564} (2003), 119--166. 
\par }
\vskip 5mm
{\sevenrm\baselineskip9pt
G\'erald Tenenbaum\par
Institut \'Elie Cartan\par 
Universit\'e de Lorraine\par
 BP 70239\par
54506 Vand\oe uvre Cedex\par
 France
\smallskip
internet: \seventt gerald.tenenbaum@univ-lorraine.fr\par}
\bye

%% file: itenn.tex

  
  %
  %

  %
  %
  \let\footnotea=\footnote
  \def\anote#1#2{\footnotea{\hbox{$^{#1}$}}{\eightpoint#2}}  
  \catcode`@=12 

 \def\defrefnote#1{\definexref{#1}{{\the\footnotenumber}}{refnotes}}

  %
  %


 \input eplain.tex
\makeatletter
\def\numberedfootnote{%
ÊÊ\global\advance\footnotenumber by 1
ÊÊ\@eplainfootnote{{\number\footnotenumber}}%
}%
\def\makecolumns#1/#2 {\par \begingroup
ÊÊ \@columndepth = #1
ÊÊ \advance\@columndepth by -1
ÊÊ \divide \@columndepth by #2
ÊÊ \advance\@columndepth by 1
ÊÊ \@linestogoincolumn = \@columndepth
ÊÊ \@linestogo = #1
ÊÊ \currentcolumn = 1
ÊÊ \def\@endcolumnactions{%
ÊÊÊÊÊÊ\ifnum \@linestogo<2
ÊÊÊÊÊÊÊÊ \the\crtok \egroup \endgroup \par 
ÊÊÊÊÊÊ\else
ÊÊÊÊÊÊÊÊ \global\advance\@linestogo by -1
ÊÊÊÊÊÊÊÊ \ifnum\@linestogoincolumn<2
ÊÊÊÊÊÊÊÊÊÊÊÊ\global\advance\currentcolumn by 1
ÊÊÊÊÊÊÊÊÊÊÊÊ\global\@linestogoincolumn = \@columndepth
ÊÊÊÊÊÊÊÊÊÊÊÊ\the\crtok
ÊÊÊÊÊÊÊÊ \else
ÊÊÊÊÊÊÊÊÊÊÊÊ&\global\advance\@linestogoincolumn by -1
ÊÊÊÊÊÊÊÊ \fi
ÊÊÊÊÊÊ\fi
ÊÊ }%
ÊÊ \makeactive\^^M
ÊÊ \letreturn \@endcolumnactions
ÊÊ \@columnwidth = \hsize
ÊÊÊÊ \advance\@columnwidth by -\parindent
ÊÊÊÊ \divide\@columnwidth by #2
ÊÊ \penalty\abovecolumnspenalty
ÊÊ \noindent 
ÊÊ \valign\bgroup
ÊÊÊÊ &\hbox to \@columnwidth{\strut \hsize = \@columnwidth ##\hfil}\cr
}%
\makeatother

\lefteqnumbers
   \def\testd{oui}
   \def\choixlat{\ifx\numadroite\testd\righteqnumbers
            \else  \lefteqnumbers\fi}
    \choixlat

\catcode`@=\letter
\def\@eplainfootnote#1{\let\@sf\empty 
  \ifhmode\edef\@sf{\spacefactor\the\spacefactor}\/\fi
  \global\advance\hlfootlabelnumber by 1
  \hlstart@impl{foot}{\hlfootlabel}%
  \hldest@impl{footback}{\hlfootbacklabel}%
  \hbox{$^{(#1)}$}%
  \hlend@impl{foot}%
  \@sf\vfootnote{#1.}%
}%
\catcode`@=\other

  \interfootnoteskip=0pt
  \let\note=\numberedfootnote
  \everyfootnote={\eightpoint\leftskip=5truemm\rightskip5truemm}
  
  \hsize150truemm\vsize 240truemm\hoffset=5truemm
  \def\dimstand{\hsize 150truemm\vsize 240truemm\hoffset=5truemm\voffset=0truemm}
  
  \pretolerance=500\tolerance=1000\brokenpenalty=5000
  \parindent3mm
  
  \countdef\temps=170
  \temps=\time
  \countdef\nminutes=171{\nminutes = \time}
  \countdef\nheure=172
  \def\heure{\begingroup                     
     \temps = \time \divide\temps by 60
     \nheure = \temps                        
     \nminutes = \time
     \multiply\temps by 60
     \advance\nminutes by -\temps            
     \ifnum\nminutes<10 \toks1 = {0}%
     \else\toks1 = {}%
     \fi
     \number\nheure h\the\toks1 \number\nminutes  
  \endgroup}%

  \newcount\chstart
  \chstart=\pageno
 \headline={\ifnum\pageno=\chstart {\hfill} \else {\hss \tenrm --\ \folio\ --\hss}\fi}
  \footline={\hfill}
  \normalbaselines
  \frenchspacing
    \def\dater{\vglue-10mm\rightline{(\the\day/\the\month/\the\year)}}
  \def\dateheure{\vglue-10mm\rightline{(\the\day/\the\month/\the\year,\ \heure)}}

  \newif\ifpagetitre \pagetitretrue
\newtoks\hautpagetitre \hautpagetitre={\hfill}
\newtoks\baspagetitre \baspagetitre={\hfill}
\newtoks\auteurcourant \auteurcourant={\hfill}
\newtoks\titrecourant \titrecourant={\hfill}
\newtoks\hautpagegauche
\newtoks\hautpagedroite
\newtoks\hautpagemilieu
\hautpagemilieu={\tenrm\hfil -- \folio\ -- \hfil}
\hautpagegauche={\ifx\midfolio\oui\the\hautpagemilieu\else\tenrm\folio\hfill\the\auteurcourant\hfill\fi}
\hautpagedroite={\ifx\midfolio\oui\the\hautpagemilieu\else\hfill\the\titrecourant\hfill\tenrm\folio\fi}
\newtoks\baspagegauche \baspagegauche={\hfil}
\newtoks\baspagedroite \baspagedroite={\hfil}
\headline={\ifpagetitre\the\hautpagetitre
\else\ifodd\pageno\the\hautpagedroite\else\the\hautpagegauche\fi\fi }
\footline={\ifpagetitre\the\baspagetitre
\else\ifodd\pageno\the\baspagedroite
\else\the\baspagegauche\fi\fi \global\pagetitrefalse}

\def\pageblanche{\vfill\eject\pagetitretrue
\null\vfill\eject
\pagetitretrue
}
\def\chgtpage{\ifodd\pageno \else
\pageblanche \fi \pagetitretrue\titreun=0\footnotenumber=0}

\def\chgtpageincrtitreun{\ifodd\pageno \else
\pageblanche \fi \pagetitretrue\footnotenumber=0}

\def\majnombres{\ifodd\pageno \else
\pageblanche \fi \pagetitretrue\hautpoly\titreun=0\footnotenumber=0}

\def\hautspages#1#2{\auteurcourant={\ninepcap#1}\titrecourant={\nineit#2}}

\ifnum\chstart=\pageno \pagetitretrue\fi
  

  \def\PAR{\par}
  
  \def\leftnote#1{\vadjust{\setbox1=\vtop{\hsize 20mm\parindent=0pt\eightpoint
  \baselineskip=9pt\rightskip=4mm plus 4mm\vskip-4.7mm#1}\hbox{\kern-2cm\smash{\box1}}}}

  
  \def\raggedcenter{\leftskip=20pt plus 10em  
       \rightskip=\leftskip 
        \parfillskip=0pt 
         \spaceskip=.3333em \xspaceskip=.5em 
          \pretolerance=9999 \tolerance=9999
           \hyphenpenalty=9999 \exhyphenpenalty=9999 }
           
  \def\titrecentre#1{{\parindent0mm\raggedcenter
       \spaceskip=.6em plus .2em minus .2em\xspaceskip=.6em plus .2em minus .2em
        \tit#1\par}}
        


  \def\oui{oui}
  
   \def\fontetitreun{\ifx\paradouze\oui\douzepts\gpdouze\twelvebf\textfont1=\twelveib\else
\quatorzepts\gpquatorze\fourteenbf\fi}

\def\fontetitreunl{\douzepts\textfont1=\twelveib\scriptfont1=\tenib\fourteenti}
 
 \def\fontetitredeux{\textfont1=\eleveni\ifx\paradouze\oui\onzepts\scriptfont1=\ninei\elevenit\else
                        \douzepts\twelveit\fi}
 
   \def\fontetitredeuxb{\ifx\paradouze\oui\onzepts\eleventi\gponze\textfont1=\elevenib\scriptfont1=\nineib
                         \else\douzepts\twelveti\scriptfont1=\twelveib\scriptfont1=\tenib\gpdouze\fi}
                         
\def\fontetitredeuxl{\onzepts\textfont1=\elevenbf\scriptfont1=\ninebf\twelvebf}
  
\def\fontetitretrois{\textfont0=\elevenrm\scriptfont0=\eightrm\textfont1=\eleveni
                      \scriptfont1=\eighti\scriptscriptfont1=\sixi\elevenit}
                      
\def\fontetitrequatre{\textfont0=\elevenrm\scriptfont0=\eightrm\textfont1=\eleveni
                      \scriptfont1=\eighti\scriptscriptfont1=\sixi\elevenrm}
  
  \newcount\titreun\titreun=0
  \newcount\titredeux\titredeux=0
  \newcount\titretrois\titretrois=0
  \newcount\titrequatre\titrequatre=0
  \newcount\enonce\enonce=0
  
  \def\incr#1{\global\advance#1 by 1 {\the #1}}
  \def\avance#1{\global\advance#1 by 1}
  \def\init#1{\global#1=0}
  
  \long\def\Indentation#1#2{\setbox10=\hbox{\fontetitreun#1}
  	                    \ifdim\wd10 < 4mm
                         \setbox10=\hbox to 4mm{\box10\hfill}
                       \else\ifdim\wd10 < 6mm
                         \setbox10=\hbox to 6mm{\box10\hfill}
  	                    \else\ifdim\wd10 < 8mm
                         \setbox10=\hbox to 8mm{\box10\hfill}
                       \else\ifdim\wd10 < 12mm
                         \setbox10=\hbox to 12mm{\box10\hfill}
                       \fi\fi\fi\fi
                       \dimen10=\hsize
                       \advance \dimen10 by -\wd10
                       \noindent \box10 %
                       \ignorespaces
                       \hbox{\vtop{\hsize=\dimen10\raggedright\noindent\fontetitreun#2}}}

  \long\def\paraun#1{\removelastskip\par\medskip\goodbreak\vskip0pt plus.01\vsize\penalty-100
                \vskip0pt plus-.01\vsize
  	              \init{\titredeux}\ifnum\optionparag=1{\init\eqnumber\init\enonce}\else{}\fi
                  \goodbreak{\fontetitreun
  	                \Indentation{\incr{\titreun}.\ }{\fontetitreun #1\par}}\nobreak\medskip}

 %
 %
 \long\def\paraunc#1{\removelastskip\par\bigskip\goodbreak\vskip0pt plus.01\vsize\penalty-100
                \vskip0pt plus-.01\vsize
  	              \init{\titredeux}
                 \ifnum\optionparag=1{\init{\eqnumber}\init\enonce}\else{}\fi
                  \goodbreak
  	                {\parindent0mm\raggedcenter\fontetitreun\incr{\titreun}.\ 
                     \fontetitreun #1\par}\nobreak\medskip}
                     
\newtoks\titreunl
\titreunl={\ifnum\titreun=1{I}\fi%
\ifnum\titreun=2{II}\fi%
\ifnum\titreun=3{III}\fi%
\ifnum\titreun=4{IV}\fi%
\ifnum\titreun=5{V}\fi%
\ifnum\titreun=6{VI}\fi%
\ifnum\titreun=7{VII}\fi%
\ifnum\titreun=8{VIII}\fi%
\ifnum\titreun=9{IX}\fi%
\ifnum\titreun=10{X}\fi%
\ifnum\titreun=11{XI}\fi%
\ifnum\titreun=12{XII}\fi%
\ifnum\titreun=13{XIII}\fi%
}
\long\def\paraunl#1{\removelastskip\par\bigskip\bigskip\goodbreak\vskip0pt plus.01\vsize\penalty-100
                \vskip0pt plus-.01\vsize
  	              \init{\titredeux}\ifnum\optionparag=1{\init\eqnumber\init\enonce}\else{}\fi
                  \goodbreak{\fontetitreunl
  	                \Indentation{\global\advance\titreun by 1{\the\titreunl}.\ }{\fontetitreunl #1\par}}\nobreak\smallskip}

  
  \long\def\paradeux#1{\init{\titretrois}\vskip0pt plus.01\vsize\penalty-10
                \vskip0pt plus-.01\vsize\ifx \elie\oui\medskip\ifnum\titredeux>0\medskip\fi\fi
                 \Indentation{\fontetitredeux\the\titreun${\cdot}$\incr{\titredeux}.}
                              {\fontetitredeux\textfont1=\eleveni#1}\nobreak\par }
  
  \long\def\paradeuxb#1{\init{\titretrois}\vskip0pt plus.001\vsize\penalty-10
                \vskip0pt plus-.01\vsize{\ifx \elie\oui\medskip\ifnum\titredeux>0\medskip\fi\fi
                  \Indentation
  {\fontetitredeuxb\the\titreun${\cdot}$\incr{\titredeux}.}{ \fontetitredeuxb#1}}\nobreak
\smallskip}

\newtoks\titredeuxl
\titredeuxl={\ifnum\titredeux=1{A}\fi%
\ifnum\titredeux=2{B}\fi%
\ifnum\titredeux=3{C}\fi%
\ifnum\titredeux=4{D}\fi%
\ifnum\titredeux=5{E}\fi%
\ifnum\titredeux=6{F}\fi%
\ifnum\titredeux=7{G}\fi%
\ifnum\titredeux=8{H}\fi%
\ifnum\titredeux=9{I}\fi%
\ifnum\titredeux=10{J}\fi%
\ifnum\titredeux=11{K}\fi%
\ifnum\titredeux=12{L}\fi%
\ifnum\titredeux=13{M}\fi%
}
 \long\def\paradeuxl#1{\init{\titretrois}\vskip0pt plus.001\vsize\penalty-10
                \vskip0pt plus-.01
                \vsize \bigskip%
                  \Indentation
     {\fontetitredeuxl\global\advance\titredeux by 1
  \quad \the\titreunl${\cdot}$\the\titredeuxl.}{ \fontetitredeuxl#1}
  \removelastskip\nobreak\smallskip}
  

  \long\def\paratrois#1{\init{\titrequatre}\ifdim\lastskip<\smallskipamount
                \removelastskip\smallskip\fi
                 \vskip0pt plus.01\vsize\penalty-10
                  \vskip0pt
plus-.01\vsize{\ifx \elie\oui\ifnum\titretrois>0\medskip\fi\fi
\Indentation{\fontetitretrois\the\titreun${\cdot}$\the\titredeux${\cdot}$\incr{\titretrois}.\ }
  {\hskip0mm\baselineskip=14pt\fontetitretrois#1}\nobreak\smallskip}}
  
  
  \long\def\paratroisl#1{\init{\titrequatre}\ifdim\lastskip<\smallskipamount
                \removelastskip\fi
                 \vskip0pt plus.01\vsize\penalty-10
                  \vskip0pt
plus-.01\vsize\ifx \elie\oui\bigskip
\fi
\Indentation{\fontetitretrois\quad \quad \the\titreunl{${\cdot}$}\the\titredeuxl${\cdot}$\incr{\titretrois}.\ }
  {\hskip0mm\fontetitretrois#1}\nobreak\smallskip}


  \long\def\paraquatre#1{\ifdim\lastskip<\smallskipamount
                \removelastskip\smallskip\fi
                 \vskip0pt plus.01\vsize\penalty-10
                  \vskip0pt
                  plus-.01\vsize\par
 
\Indentation{\fontetitrequatre \the\titreun{${\cdot}$}\the\titredeux${\cdot}$\the\titretrois${\cdot}$\incr{\titrequatre}.\ }
{\hskip0mm\fontetitrequatre#1}\nobreak\smallskip}


\newtoks\titrequatrel
\titrequatrel={\ifnum\titrequatre=1{a}\fi%
\ifnum\titrequatre=2{b}\fi%
\ifnum\titrequatre=3{c}\fi%
\ifnum\titrequatre=4{d}\fi%
\ifnum\titrequatre=5{e}\fi%
\ifnum\titrequatre=6{f}\fi%
\ifnum\titrequatre=7{g}\fi%
\ifnum\titrequatre=8{h}\fi%
\ifnum\titrequatre=9{i}\fi%
\ifnum\titrequatre=10{j}\fi%
\ifnum\titrequatre=11{k}\fi%
\ifnum\titrequatre=12{l}\fi%
\ifnum\titrequatre=13{m}\fi%
}
\long\def\paraquatrel#1{\ifdim\lastskip<\smallskipamount
                \removelastskip\smallskip\fi
                 \vskip0pt plus.01\vsize\penalty-10
                  \vskip0pt
                  plus-.01\vsize{\bigskip
\Indentation{\global\advance\titrequatre by 1
\fontetitrequatre\quad \quad \quad \the\titreunl${\cdot}$\the\titredeuxl${\cdot}$\the\titretrois${\cdot}$\the\titrequatrel.\ }
{\hskip0mm\fontetitrequatre#1}\nobreak\smallskip}}

\ifx\optionkeys\oui
\def\drefun#1{\definexref{¤#1}{{\the\titreun}}{}} 
\def\drefdeux#1{\definexref{¤#1}{{\the\titreun}.{\the\titredeux}}{}}
\def\dreftrois#1{\definexref{¤#1}{{\the\titreun}.{\the\titredeux}.{\the\titretrois}}{}}
\else
\def\drefun#1{\definexref{prg#1}{{\the\titreun}}{}} 
\def\drefdeux#1{\definexref{prg#1}{{\the\titreun}.{\the\titredeux}}{}}
\def\dreftrois#1{\definexref{prg#1}{{\the\titreun}.{\the\titredeux}.{\the\titretrois}}{}}
\fi

%


  \long\def\propdeux#1#2#3#4{%
       \avance{\enonce}
       \leavevmode\edef\temp{#2}%
         \ifx\temp\empty 
          \else
           \definexref{#2}{#1~{\the\titreun.\the\enonce}}{enonces}
            \definexref{s#2}{{\the\titreun.\the\enonce}}{enonces}
             \fi
\smallskip
      \noindent{\bf#1\ {\bf\the\titreun.\the\enonce{#3}.}\enspace}{\sl#4\par}%
      \ifdim\lastskip<\medskipamount \removelastskip\penalty55\par \fi
   }

  \long\def\propun#1#2#3#4{%
      \avance{\enonce}
       \leavevmode\edef\temp{#2}%
        \ifx\temp\empty 
          \else
           \definexref{#2}{#1~{\the\enonce}}{enonces}
            \definexref{{s#2}}{{\the\enonce}}{enonces}
             \fi
   \par 
     \noindent{\bf#1\ {\bf\the\enonce{#3}.}\enspace}{\sl#4\par}%
     \ifdim\lastskip<\medskipamount \removelastskip\penalty55\medskip\fi
  }
  
  \long\def\prop#1#2#3#4{\ifnum\optionparag=1
                          \propdeux{#1}{#2}{\textfont1=\elevenib#3}{#4} \else\propun{#1}{#2}{\textfont1=\elevenib#3}{#4}\fi}

  \long\def\propt#1#2#3{\ifx\tpf\oui \prop{Th\'eo\-r\`eme}{#1}{#2}{#3}\par
                       \else\prop{Theorem}{#1}{#2}{#3}\par\fi}
  \long\def\Propt#1#2{\propt{#1}{}{#2}}
  \long\def\propl#1#2#3{\ifx\tpf\oui\prop{Lem\-me}{#1}{#2}{#3}\par
                         \else\prop{Lemma}{#1}{#2}{#3}\par\fi}
  
  \long\def\propc#1#2#3{\ifx\tpf\oui\prop{Corol\-laire}{#1}{#2}{#3}\par
                         \else\prop{Corollary}{#1}{#2}{#3}\par\fi}

  \long\def\propd#1#2#3{\ifx\tpf\oui\prop{D\'efi\-nition}{#1}{#2}{#3}\par
                       \else\prop{Definition}{#1}{#2}{#3}\par\fi} 
  
  \long\def\proptd#1#2#3{\ifx\tpf\oui\prop{Th\'eor\`eme et d\'efi\-nition}{#1}{#2}{#3}\par
                       \else\prop{Theorem and definition}{#1}{#2}{#3}\par\fi}


  
  \newcount\optionparag\optionparag=1
  
  \long\def\section#1#2{\ifnum\optionparag=1 \paraun{#2} 
                        \else\goodbreak{\fontetitreun
  	                \Indentation{#1.\ }{#2}}\nobreak\smallskip\fi}

  \def\eqconstruct#1{\ifnum\optionparag=1{\the\titreun\hbox{$\cdot$}#1}\else{#1}\fi}

  
  
  \def\numref{oui}  
  
  \newcount\mesref\mesref=0 
  \def\defbib#1{\ifx\numref\oui\global\advance\mesref by 1 \definexref{#1}{{\the
                 \mesref}}{}\else\definexref{#1}{#1}{}\fi}
  \def\bibtem#1{\defbib{#1}\item{\citer{#1}}}
  \def\citer#1{[\ref{#1}]}

  
  \font\seventeenmsa=msam10 at 17pt    
  \font\fourteenmsa=msam10 at 14pt
  \font\twelvemsa=msam10 at 12pt
  \font\tenmsa=msam10                 
  \font\ninemsa=msam10 at 9pt 
  \font\eightmsa=msam10 at 8pt 
  \font\sevenmsa=msam7 
  \font\sixmsa=msam10 at 6pt
  \font\fivemsa=msam5
  \newfam\msafam\textfont\msafam=\tenmsa\scriptfont\msafam=\sevenmsa\scriptscriptfont\msafam=\fivemsa
  
  \font\seventeenbb=msbm10 at 17pt     
  \font\fourteenbb=msbm10 at 14pt
  \font\twelvebb=msbm10 at 12pt
  \font\tenbb=msbm10                   
  \font\ninebb=msbm10 at 9pt 
  \font\eightbb=msbm10 at 8pt 
  \font\sevenbb=msbm7 
  \font\sixbb=msbm10 at 6pt
  \font\fivebb=msbm5 
  \newfam\bbfam\textfont\bbfam=\tenbb\scriptfont\bbfam=\sevenbb\scriptscriptfont\bbfam=\fivebb
  \def\bb{\fam\bbfam\tenbb}%

  \font\seventeenscaln=eusm10 at 17pt   
  \font\twelvescaln=eusm10 at 12pt
  \font\tenscaln=eusm10                
  \font\ninescaln=eusm10 scaled 900
  \font\eightscaln=eusm10 scaled 800
  \font\sevenscaln=eusm10 scaled 700
  \font\sixscaln=eusm10 scaled 600
   
  \newfam\scalnfam\textfont\scalnfam=\tenscaln\scriptfont\scalnfam=\sevenscaln\scriptscriptfont\scalnfam=\sixscaln
  \def\scaln{\fam\scalnfam\tenscaln}%
  \def\scal{\scaln}
  
  \font\tenscalb=eusb10                

  \font\sevenscalb=eusb10 scaled 700

  \newfam\scalbfam\textfont\scalbfam=\tenscalb\scriptfont\scalbfam=\sevenscalb
  %
  
  %
  %
  \font\fourteenrm=cmr12 scaled 1200
  \font\elevenrm=cmr10 at 11pt
  \font\twelverm=cmr12
  \font\ninerm=cmr9
  \font\eightrm=cmr8      
  \font\sevenrm=cmr7
  \font\sixrm=cmr6

  \font\seventeenpcap=cmcsc10 at 17pt
  \font\tenpcap=cmcsc10                        
  \font\ninepcap=cmcsc9
  \font\eightpcap=cmcsc8
  \font\sevenpcap=cmcsc10 scaled 700
  
  \newfam\pcapfam\textfont\pcapfam=\tenpcap\scriptfont\pcapfam=\sevenpcap
  \def\pcap{\fam\pcapfam\tenpcap}
  
  \font\seventeenrm=cmbx12 scaled 1400

  \font\fourteenbf=cmbx10 scaled 1400
  
  \font\twelvebf=cmbx12
  \font\elevenbf=cmbx10 at 11pt
  \font\ninebf=cmbx9  
  \font\eightbf=cmbx8
  \font\sixbf=cmbx6
  
  \font\tengot=eufm10                           
   
  \font\eightgot=eufm10 at 8truept 
  \font\sevengot=eufm7 
  \font\sixgot=eufm10 at 6 truept 
   
  \newfam\gotfam
  \textfont\gotfam=\tengot\scriptfont\gotfam=\sevengot\scriptscriptfont\gotfam=\sixgot
  %

  
  \def\tit{%
  \textfont0=\seventeenrm\scriptfont0=\tenrm\def\rm{\fam0\seventeenrm}%
  \textfont1=\seventeenib\scriptfont1=\twelveib%
  \textfont2=\seventeensy\scriptfont2=\twelvesy\scriptscriptfont2=\ninesy
  \textfont3=\seventeenex
  \textfont\itfam=\seventeenti
  \def\it{\fam\itfam\seventeenti}%
  \textfont\bbfam=\seventeenbb \scriptfont\bbfam=\twelvebb
  \def\bb{\fam\bbfam\seventeenbb}%
  \textfont\msafam=\seventeenmsa\scriptfont\msafam=\twelvemsa
  \textfont\scalnfam=\seventeenscaln
  \def\pcap{\fam\pcapfam\seventeenpcap}
  \normalbaselineskip=25pt\normalbaselines\rm}

  \font\seventeenti=cmbxti10 scaled 1680
  
  \font\fourteenti=cmbxti10 at 14pt
  
  \font\twelveti=cmbxti10 scaled 1200
  \font\eleventi=cmbxti10 at 11pt

  %
  %
  \font\twelveit=cmti12	
  \font\elevenit=cmti10 scaled 1100
  \font\nineit=cmti9
  \font\eightit=cmti8
  \font\sevenit=cmti7

  %
  %
 
 \font\seventeenib=cmmib10 scaled 1680
  \font\fourteenib=cmmib10 scaled 1400
  \font\twelveib=cmmib10 scaled 1200
  \font\elevenib=cmmib10 scaled 1100
  \font\tenib=cmmib10
\font\eightib=cmmib10 scaled 800
  \font\nineib=cmmib10 scaled 900
\font\sevenib=cmmib10 scaled 700
\font\sixib=cmmib10 scaled 600
\font\fiveib=cmmib10 scaled 500

\ifx\ITAN\oui
\else
\innernewfam\cmmibfam
\textfont\cmmibfam=\tenib
\scriptfont\cmmibfam=\sevenib
\scriptscriptfont\cmmibfam=\fiveib
\def\ib{\fam\cmmibfam\tenib}
\fi

  %
  %
  
  \font\eleveni=cmmi10 scaled 1100
  \font\ninei=cmmi9
  \font\eighti=cmmi8 
  \font\seveni=cmmi7 			                
  \font\sixi=cmmi6
  
  \font\ninesl=cmsl9                    
  \font\eightsl=cmsl8 
  \font\sevensl=cmsl10 at 7pt

  \font\ninett=cmtt9                    
  \font\eighttt=cmtt8
  \font\seventt=cmtt10 scaled 700

  \font\seventeensy=cmsy10 scaled 1680    
  \font\fourteensy=cmsy10 scaled 1400
  \font\twelvesy=cmsy10 scaled 1176
  
  \font\ninesy=cmsy9                      
  \font\eightsy=cmsy8
  \font\sixsy=cmsy6
  \font\seventeenex=cmex10 at 17pt
  \font\fourteenex=cmex10 at 14pt
  \font\twelveex=cmex10 at 12pt
  \font\nineex=cmex10 at 9pt
  \font\eightex=cmex10 at 8pt
  \font\sevenex=cmex10 at 7pt
  \font\sixex=cmex10 at 6pt
  \font\fiveex=cmex10 at 5pt
  
   
  \font\fourteengp=cmmi10 at 14pt
  
  \font\twelvegp=cmmib10 at 12pt
  \font\elevengp=cmmib10 at 11pt
  \font\tengp=cmmib10                          
  \font\ninegp=cmmib10 at 9pt 
  \font\eightgp=cmmib8 
   
  \font\sixgp=cmmib6


  \def\gponze{\textfont0=\elevenbf\scriptfont0=\eightbf\scriptscriptfont0=\sixbf
           \textfont1=\elevengp\scriptfont1=\eightgp\scriptscriptfont1=\sixgp}
  \def\gpdouze{\textfont0=\twelvebf\scriptfont0=\tenbf\scriptscriptfont0=\ninebf
           \textfont1=\twelvegp\scriptfont1=\tengp\scriptscriptfont1=\ninegp}        
  
 \def\gpquatorze{\textfont0=\fourteenbf\scriptfont0=\twelvebf\scriptscriptfont0=\elevenbf
           \textfont1=\fourteengp\scriptfont1=\twelvegp\scriptscriptfont1=\elevengp}

  
  \expandafter\chardef\csname pre amssym.def at\endcsname=\the\catcode`\@
  \catcode`\@=11
  \def\undefine#1{\let#1\undefined}
  \def\newsymbol#1#2#3#4#5{\let\next@\relax
   \ifnum#2=\@ne\let\next@\msafam@\else
   \ifnum#2=\tw@\let\next@\bbfam@\fi\fi
   \mathchardef#1="#3\next@#4#5}
  \def\mathhexbox@#1#2#3{\relax
   \ifmmode\mathpalette{}{\m@th\mathchar"#1#2#3}%
   \else\leavevmode\hbox{$\m@th\mathchar"#1#2#3$}\fi}
  \def\hexnumber@#1{\ifcase#1 0\or 1\or 2\or 3\or 4\or 5\or 6\or 7\or 8\or
   9\or A\or B\or C\or D\or E\or F\fi}
  
  \def\setboxz@h{\setbox\z@\hbox}
  \def\wdz@{\wd\z@}
  \def\boxz@{\box\z@}
  
  \edef\msafam@{\hexnumber@\msafam}
  \mathchardef\dabar@"0\msafam@39
  
  \edef\bbfam@{\hexnumber@\bbfam}
  \def\widehat#1{\setboxz@h{$\m@th#1$}%
   \ifdim\wdz@>\tw@ em\mathaccent"0\bbfam@5B{#1}%
   \else\mathaccent"0362{#1}\fi}
  \def\widetilde#1{\setboxz@h{$\m@th#1$}%
   \ifdim\wdz@>\tw@ em\mathaccent"0\bbfam@5D{#1}%
   \else\mathaccent"0365{#1}\fi}
  \newsymbol\leqq 1335          
  \newsymbol\leqslant 1336
  \newsymbol\lessgtr 1337       
  \newsymbol\backprime 1038     
  \newsymbol\risingdotseq 133A  
  \newsymbol\fallingdotseq 133B 
  \newsymbol\succcurlyeq 133C   
  \newsymbol\geqq 133D          
  \newsymbol\geqslant 133E
  \newsymbol\nmid 232D
  \newsymbol\nexists 2040
  \newsymbol\smallsetminus 2272
  \newsymbol\varnothing 203F
  
  \catcode`\@=\active

  \catcode`\@=11
  \newcount\typofr\typofr=1
  
  \catcode`\;=\active
  \def;{\ifnum\typofr=1\relax\ifhmode\ifdim\lastskip>\z@\unskip\fi
     \kern.2em\fi\string;\else\string;\fi}
  
  \catcode`\:=\active
  \def:{\ifnum\typofr=1\relax\ifhmode\ifdim\lastskip>\z@\unskip\fi
  \penalty\@M\ \fi\string:\else\string:\fi}
  
  \catcode`\!=\active
  \def!{\ifnum\typofr=1\relax\ifhmode\ifdim\lastskip>\z@\unskip\fi
     \kern.2em\fi\string!\else\string!\fi}
  
  \catcode`\?=\active
  \def?{\ifnum\typofr=1\relax\ifhmode\ifdim\lastskip>\z@\unskip\fi
     \kern.2em\fi\string?\else\string?\fi}

  \def\francais{\typofr=1\def\tpf{oui}}
  \def\anglais{\typofr=2\def\tpf{non}\def\english{oui}}
  \def\oui{oui}
  \francais
  
  \catcode`\@=12
  

\ifx\textures\oui
\def\raye #1|{\leavevmode\setbox1=\hbox{#1}%
\raise .5pt\hbox to \wd1{\xleaders\hbox{\rge{$ \sct / $}%
\kern 1pt}\hfill\kern -1pt }\kern -\wd1 \unhbox1\relax }

\def\barre #1|{\leavevmode\setbox1=\hbox{#1}%
\rlap{\Red\vrule height 2.4pt depth -1.2pt width \wd1}\Black \unhbox1\relax}
\else
\def\raye #1|{\leavevmode\setbox1=\hbox{#1}%
\raise .5pt\hbox to \wd1{\xleaders\hbox{\rge{$ \sct / $}%
\kern 1pt}\hfill\kern -1pt }\kern -\wd1 \unhbox1\relax }

\def\barre #1|{\leavevmode\setbox1=\hbox{#1}%
\rlap{\color{red}\vrule height 2.4pt depth -1.2pt width \wd1}\color{black} \unhbox1\relax}

\fi
  

  
  \def\og{\leavevmode\raise.24ex\hbox{$\scriptscriptstyle\langle\!\langle\>$}}    
  \def\fg{\leavevmode\raise.24ex\hbox{$\scriptscriptstyle\>\rangle\!\rangle$}}    

  \def\d{\,{\rm d}}

  \def\CC{{\bb C}}

  \def\O{{\scal O}}
  \def\P{{\scaln P}}

  \def\frac#1#2{{#1\over #2}}
  \def\di#1#2{\sct#1\atop{\sct#2}}
  \def\tri#1#2#3{{\sct#1\atop\sct#2}\atop\sct#3}

  \def\numero{n$^{\rm o}\thinspace$}
\def\numeros{n$^{\rm os}\thinspace$}


  \def\numero{n$^{\rm o}\thinspace$}

  \def\¤{\S\thinspace}

  \def\¥{$\bullet$ }
  
  
  \def\e{{\rm e}}

  \def\epsilon{\varepsilon}

  \def\phi{\varphi}
  \def\theta{\vartheta}
  \def\rho{\varrho}
  \def\dm{{\textstyle{1\over 2}}}
  \def\txt{\textstyle}
  
  \def\sct{\scriptstyle}
  \def\pf{\noi{\it Proof. }}
  \def\nid{\ifnum\typofr=1\par\noindent{\it D\'emonstration. }\else\pf\fi}
  \def\noi{\noindent}
  \def\rem{\ifnum\typofr=1\noi{\it Remarque.}\ \else\noi{\it Remark.}\ \fi}
  \def\rems{\ifnum\typofr=1\noi{\it Remarques.}\ \else\noi{\it Remarks.}\ \fi}

  \def\1{{\bf 1}}
  \def\|{\Vert}

  \def\leq{\leqslant}
  \def\geq{\geqslant}

  \def\ie{{i.e.\ }}
  \def\eg{{e.g.}}
  

  \def\log{\mathop{\rm log}\nolimits}
  \def\ft#1#2{{\txt{#1\over #2}}}




  \def\pmb#1{\setbox0=\hbox{#1}%
  \kern-.025em\copy0\kern-\wd0\kern.05em\copy0\kern-\wd0\kern-.025em\raise .0433em\box0 }

  
  \skewchar\eighti='177 \skewchar\sixi='177
  \skewchar\eightsy='60 \skewchar\sixsy='60
  
  \def\eightpoint{%
  \textfont0=\eightrm\scriptfont0=\sixrm\scriptscriptfont0=\fiverm
  \def\rm{\fam0\eightrm}%
  \textfont1=\eighti\scriptfont1=\sixi
  \scriptscriptfont1=\fivei\def\oldstyle{\fam1\seveni}%
  \textfont2=\eightsy\scriptfont2=\sixsy\scriptscriptfont2=\fivesy
  \textfont3=\eightex\scriptfont3=\sixex
  \textfont\itfam=\eightit
  \def\it{\fam\itfam\eightit}%
  \textfont\slfam=\eightsl
  \def\sl{\fam\slfam\eightsl}%
  \textfont\bbfam=\eightbb \scriptfont\bbfam=\sixbb\scriptscriptfont\bbfam=\fivebb
  \def\bb{\fam\bbfam\eightbb}%
  \textfont\msafam=\eightmsa\scriptfont\msafam=\sixmsa
  \textfont\scalnfam=\eightscaln
  \def\scaln{\fam\scalnfam\eightscaln}
  \textfont\ttfam=\eighttt
  \def\tt{\fam\ttfam\eighttt}%
\textfont\gotfam=\eightgot
  \textfont\bffam=\eightbf\scriptfont\bffam=\sixbf\scriptscriptfont\bffam=\fivebf
  \def\bf{\fam\bffam\eightbf}%
  \ifx\ITAN\oui\else\textfont\cmmibfam=\eightib
       \scriptfont\cmmibfam=\sixib
        \scriptscriptfont\cmmibfam=\fiveib
         \def\ib{\fam\cmmibfam\eightib}
   \fi
  \textfont\pcapfam=\eightpcap
  \def\pcap{\fam\pcapfam\eightpcap}
  \abovedisplayskip=2pt plus2pt minus 2pt
  \belowdisplayskip=2pt plus1pt minus 2pt
  \abovedisplayshortskip= 1pt plus 2pt minus 1pt
  \belowdisplayshortskip= 1pt plus 2pt minus 1pt
  \smallskipamount=2pt plus 1pt minus 2pt
  \medskipamount=3pt plus 2pt minus 2pt
  \bigskipamount=7pt plus 3pt minus 3pt
  \setbox\strutbox=\hbox{\vrule height 5pt depth 2pt width 0pt}%
  \normalbaselineskip=9pt\normalbaselines\rm}

  \def\({\left(}
  \def\){\right)}
  
  \def\footnoterule{\kern -2pt\hrule width 7truecm\kern 2.4pt}
  
  \def\xnotedef#1{\definexref{#1}{\noexpand\number\footnotenumber}{Note}}%

  
  
  \def\ninepoint{%
  \textfont0=\ninerm\scriptfont0=\sixrm\scriptscriptfont0=\fiverm
  \def\rm{\fam0\ninerm}%
  \textfont1=\ninei\scriptfont1=\sixi
  \scriptscriptfont1=\fivei\def\oldstyle{\fam1\ninei}%
  \textfont2=\ninesy\scriptfont2=\sixsy\scriptscriptfont2=\fivesy
  \textfont3=\nineex\scriptfont3=\sixex
  \textfont\itfam=\nineit
  \def\it{\fam\itfam\nineit}%
  \textfont\slfam=\ninesl
  \def\sl{\fam\slfam\ninesl}%
  \textfont\bbfam=\ninebb\scriptfont\bbfam=\sixbb\scriptscriptfont\bbfam=\fivebb
  \def\bb{\fam\bbfam\ninebb}%
  \textfont\msafam=\ninemsa\scriptfont\msafam=\sixmsa\scriptscriptfont\msafam=\fivemsa
  \textfont\scalnfam=\ninescaln
  \def\scaln{\fam\scalnfam\ninescaln}
  \textfont\ttfam=\ninett
  \def\tt{\fam\ttfam\ninett}%
  \textfont\bffam=\ninebf\scriptfont\bffam=\sixbf\scriptscriptfont\bffam=\fivebf
  \def\bf{\fam\bffam\ninebf}%
  \abovedisplayskip=3pt plus2pt minus 2pt
  \belowdisplayskip=3pt plus1pt minus 2pt
  \abovedisplayshortskip= 2pt plus 2pt minus 1pt
  \belowdisplayshortskip= 2pt plus 2pt minus 1pt
  \smallskipamount=2pt plus 1pt minus 2pt
  \medskipamount=3pt plus 2pt minus 2pt
  \bigskipamount=7pt plus 3pt minus 3pt
  \setbox\strutbox=\hbox{\vrule height 5pt depth 2pt width 0pt}%
  \normalbaselineskip=10.5pt plus.3pt minus.3pt\normalbaselines\rm}

  \def\sevenpoint{%
  \textfont0=\sevenrm\scriptfont0=\sixrm\scriptscriptfont0=\fiverm
  \def\rm{\fam0\sevenrm}%
  \textfont1=\seveni\scriptfont1=\sixi
  \scriptscriptfont1=\fivei\def\oldstyle{\fam1\seveni}%
  \textfont2=\sevensy\scriptfont2=\sixsy\scriptscriptfont2=\fivesy
  \textfont3=\sevenex\scriptfont3=\fiveex
  \textfont\itfam=\sevenit
  \def\it{\fam\itfam\sevenit}%
  \textfont\slfam=\sevensl
  \def\sl{\fam\slfam\sevensl}%
  \textfont\bbfam=\sevenbb \scriptfont\bbfam=\sixbb\scriptscriptfont\bbfam=\fivebb
  \def\bb{\fam\bbfam\sevenbb}%
  \textfont\msafam=\sevenmsa\scriptfont\msafam=\sixmsa
  \textfont\scalnfam=\sevenscaln
  \def\scaln{\fam\scalnfam\sevenscaln}
  \textfont\bffam=\sevenbf\scriptfont\bffam=\sixbf\scriptscriptfont\bffam=\fivebf
  \def\bf{\fam\bffam\sevenbf}%
  \textfont\ttfam=\seventt
  \abovedisplayskip=2pt plus2pt minus 2pt
  \belowdisplayskip=2pt plus1pt minus 2pt
  \abovedisplayshortskip= 1pt plus 2pt minus 1pt
  \belowdisplayshortskip= 1pt plus 2pt minus 1pt
  \smallskipamount=2pt plus 1pt minus 2pt
  \medskipamount=3pt plus 2pt minus 2pt
  \bigskipamount=7pt plus 3pt minus 3pt
  \setbox\strutbox=\hbox{\vrule height 5pt depth 2pt width 0pt}%
  \normalbaselineskip=9pt\normalbaselines\rm}

 \def\onzepts{%
 \textfont0=\elevenrm\scriptfont0=\ninerm
 \textfont1=\elevenib\scriptfont1=\ninei}

\def\douzepts{%
  \textfont0=\twelverm\scriptfont0=\tenrm\def\rm{\fam0\twelverm}%
  \textfont1=\twelveib\scriptfont1=\teni%
  \textfont2=\twelvesy\scriptfont2=\tensy\scriptscriptfont2=\eightsy
  \textfont3=\twelveex
  \textfont\itfam=\twelveti
  \def\it{\fam\itfam\twelveti}%
  \textfont\bffam=\twelvebf\scriptfont\bffam=\tenbf\scriptscriptfont\bffam=\eightbf
  \def\bf{\fam\bffam\twelvebf}%
  \textfont\bbfam=\twelvebb \scriptfont\bbfam=\tenbb
  \def\bb{\fam\bbfam\twelvebb}%
  \textfont\msafam=\twelvemsa\scriptfont\msafam=\tenmsa
  \textfont\scalnfam=\twelvescaln
  \normalbaselineskip=15pt\normalbaselines\rm}

\def\quatorzepts{%
  \textfont0=\fourteenrm\scriptfont0=\twelverm\def\rm{\fam0\fourteenrm}%
  \textfont1=\fourteenib\scriptfont1=\twelveib%
  \textfont2=\fourteensy\scriptfont2=\twelvesy\scriptscriptfont2=\tensy
  \textfont3=\fourteenex
  \textfont\itfam=\fourteenti
  \def\it{\fam\itfam\fourteenti}%
  \textfont\bffam=\fourteenbf\scriptfont\bffam=\twelvebf\scriptscriptfont\bffam=\tenbf
  \def\bf{\fam\bffam\fourteenbf}%
  \textfont\bbfam=\fourteenbb \scriptfont\bbfam=\twelvebb
  \def\bb{\fam\bbfam\fourteenbb}%
  \textfont\msafam=\fourteenmsa\scriptfont\msafam=\twelvemsa
  \textfont\scalnfam=\twelvescaln
  \normalbaselineskip=18pt\normalbaselines\rm}


\def\AA{{\it Acta Arith.}}

\def\picture #1 by #2 (#3){\leavevmode\vbox to #2{
     \hrule width #1 height 0pt depth 0pt
      \vfill
       \special{picture #3}}}

\def\illustration #1 by #2 (#3) scaled #4{\dimen1=#2
  \divide\dimen1 by 1000
  \multiply\dimen1 by #4
  \vtop to \dimen1{\dimen1=#1
  \divide\dimen1 by 1000
  \multiply\dimen1 by #4
  \hsize=\dimen1\vss
  \noindent\special{illustration #3 scaled #4}}}

%% file: option_keys.tex
\ifx\optionkeymacros\undefined\else \fi

\catcode`\Œ=\active\defŒ{{\aa}}       
\catcode`\º=\active\defº{\int}        
\catcode`\=\active\def{\c c}        
\catcode`\¶=\active\def¶{\partial}    
\catcode`\Ä=\active\defÄ{\oint}       
\catcode`\Æ=\active\defÆ{\triangle}   
\catcode`\Â=\active\defÂ{\neg}        
\catcode`\µ=\active\defµ{\mu}         
\catcode`\¿=\active\def¿{{\o}}        
\catcode`\¹=\active\def¹{\pi}         
\catcode`\Ï=\active\defÏ{{\oe}}       
\catcode`\§=\active\def§{{\ss}}       
\catcode`\ =\active\def {\dagger}     
\catcode`\Ã=\active\defÃ{\sqrt}       
\catcode`\·=\active\def·{\Sigma}      
\catcode`\Å=\active\defÅ{\approx}     
\catcode`\½=\active\def½{\Omega}      
\catcode`\£=\active\def£{{\it\$}}     
\catcode`\°=\active\def°{\infty}      
\catcode`\¤=\active\def¤{{\S}}        
\catcode`\¦=\active\def¦{{\P}}        
\catcode`\¥=\active\def¥{\bullet}     
\catcode`\»=\active\def»{\leavevmode\raise.585ex\hbox{\b a}}      
\catcode`\¼=\active\def¼{\leavevmode\raise.6ex\hbox{\b o}}        
\catcode`\­=\active\def­{\not=}       
\catcode`\²=\active\def²{\leq}        
\catcode`\³=\active\def³{\geq}        
\catcode`\Ö=\active\defÖ{\div}        
\catcode`\É=\active\defÉ{{\dots}}     
\catcode`\¾=\active\def¾{{\ae}}       
\catcode`\Ç=\active\defÇ{\og}         
\catcode`\Ò=\active\defÒ{``}          
\catcode`\Á=\active\defÁ{!`}          
\catcode`\¢=\active\def¢{\rlap/c}     
\catcode`\Ô=\active\defÔ{`}           
\catcode`\Õ=\active\defÕ{'}           


\catcode`\=\active\def{{\AA}}       
\catcode`\'=\active\def'{\c C}        
\catcode`\¯=\active\def¯{{\O}}        
\catcode`\¸=\active\def¸{\Pi}         
\catcode`\Î=\active\defÎ{{\OE}}       
\catcode`\®=\active\def®{{\AE}}       
\catcode`\×=\active\def×{\diamond}    
\catcode`\¡=\active\def¡{\accent'27}  
\catcode`\Ó=\active\defÓ{''}          
\catcode`\±=\active\def±{\pm}         
\catcode`\È=\active\defÈ{\fg}         
\catcode`\À=\active\defÀ{?`}          
\catcode`\Ð=\active\defÐ{--}          
\catcode`\Ñ=\active\defÑ{---}         


\catcode`\Š=\active\defŠ{\"a}        
\catcode`\'=\active\def'{\"e}        
\catcode`\•=\active\def•{\"{\i}}     
\catcode`\š=\active\defš{\"o}        
\catcode`\Ÿ=\active\defŸ{\"u}        
\catcode`\Ø=\active\defØ{\"y}        
\catcode`\å=\active\defå{\^A}        
\catcode`\€=\active\def€{\"A}        
\catcode`\…=\active\def…{\"O}        
\catcode`\†=\active\def†{\"U}        
\catcode`\‡=\active\def‡{\'a}        
\catcode`\Ž=\active\defŽ{\'e}        
\catcode`\'=\active\def'{\'{\i}}     
\catcode`\—=\active\def—{\'o}        
\catcode`\œ=\active\defœ{\'u}        
\catcode`\ƒ=\active\defƒ{\'E}        
\catcode`\æ=\active\defæ{\^E}        
\catcode`\é=\active\defé{\`E}        %
\catcode`\ˆ=\active\defˆ{\`a}        
\catcode`\=\active\def{\`e}        
\catcode`\"=\active\def"{\`{\i}}     
\catcode`\˜=\active\def˜{\`o}        
\catcode`\=\active\def{\`u}        
\catcode`\Ë=\active\defË{\`A}        
\catcode`\‹=\active\def‹{\~a}        
\catcode`\–=\active\def–{\~n}        
\catcode`\›=\active\def›{\~o}        
\catcode`\Ì=\active\defÌ{\~A}        
\catcode`\"=\active\def"{\~N}        
\catcode`\Í=\active\defÍ{\~O}        
\catcode`\‰=\active\def‰{\^a}        
\catcode`\=\active\def{\^e}        
\catcode`\"=\active\def"{\^{\i}}     
\catcode`\™=\active\def™{\^o}        
\catcode`\ž=\active\defž{\^u}        

\let\optionkeymacros\null